\documentclass{article}
\usepackage{amsfonts,amssymb}

\newcommand{\integers}{{\mathbb Z}}
\newcommand{\realnos}{{\mathbb R}}

\def\ov{\overline}

\title{Some examples of aspherical 4-manifolds that are homology 4-spheres}

\author{John G. Ratcliffe and Steven T. Tschantz}

\begin{document}





\date{Department of Mathematics, Vanderbilt University, \\  Nashville, TN 37240,  USA}


\maketitle

\begin{abstract}
In this paper, Problem 4.17 on R. Kirby's problem list is solved   
by constructing infinitely many aspherical 4-manifolds that are homology 4-spheres. 
\end{abstract}

\section{Introduction} 

Problem 4.17, attributed to W. Thurston, on R. Kirby's 1978 problem list \cite{Kirby} asks: 
\begin{quote}
Can a homology 4-sphere ever be a $K(\pi,1)$?  Who knows an example of a rational 
homology 4-sphere which is a $K(\pi,1)$?
\end{quote}
Examples of aspherical rational homology 
4-spheres were constructed by F. Luo \cite{Luo} in 1988. In this paper we affirmatively  
answer the main question of Problem 4.17 by constructing  
infinitely many aspherical 4-manifolds that are homology 4-spheres. 

In \S 1, we construct our examples by Dehn surgery on the complement of 
five linked 2-tori in the 4-sphere recently constructed by D. Ivan\v{s}i\'c \cite{Ivansic}.  
In \S 2, we reprove the main results of \S 1 by a direct homology calculation 
which is independent of Ivan\v{s}i\'c's paper \cite{Ivansic}.  
In \S 3, we derive geometric information that determines precisely 
which Dehn surgeries on Ivan\v{s}i\'c's link complement yield aspherical homology 4-spheres. 
In \S 4, we describe some nice properties of our examples. 
In particular, we discuss how recent work  
of M. Anderson \cite{Anderson} implies that many of our examples are Einstein 4-manifolds.

\section{Dehn surgery on Ivan\v{s}i\'c's linked 2-tori in the 4-sphere} 

Using J{\o}rgensen-Thurston's hyperbolic Dehn surgery theory \cite{Thurston}, \cite{Benedetti}, 
it is easy to construct infinitely many aspherical 3-manifolds that are homology 3-spheres. 
Simply perform Dehn surgery on the complement $M$ of a hyperbolic knot in $S^3$ with 
surgery coefficient of the form $1/k$, with $k$ an integer. 
The 3-manifold $M(k)$ obtained by Dehn surgery 
is a homology 3-sphere, and for all but finitely many $k$, 
the 3-manifold $M(k)$ supports a hyperbolic metric and therefore is aspherical. 
The volume of the hyperbolic 3-manifold $M(k)$ is less than 
the volume of the hyperbolic manifold $M$, and 
the volume of $M(k)$ converges to the volume of $M$ as $k$ goes to infinity. 
Hence there are infinitely many homotopy types of 3-manifolds of the form $M(k)$ 
by Mostow's rigidity theorem. 
Thus there are infinitely many aspherical 3-manifolds that are homology 3-sphere.

We construct our examples by performing Dehn surgery on the hyperbolic 
complement $M$ of five linked 2-tori in $S^4$ recently constructed 
by D. Ivan\v{s}i\'c \cite{Ivansic}. 
The hyperbolic 4-manifold $M$ is the orientable double cover of the most symmetric 
hyperbolic 4-manifold $N$ on our list \cite{Ratcliffe} of hyperbolic 4-manifolds of minimum volume. 
The Euler characteristic of $N$ is one, and so the Euler characteristic of $M$ is two. 
The hyperbolic 4-manifold $M$ has five cusps 
each of which is homeomorphic to $T^3 \times [0,\infty)$ where $T^3$ is the 3-torus.     

Here is an outline of our argument. 
We show that for infinitely many Dehn surgeries on each cusp of $M$ the closed 4-manifold 
$\hat M$ obtained by Dehn surgery is a homology 4-sphere. 
By the Gromov-Thurston $2\pi$ theorem,    
for all but finitely many surgeries on each cusp of $M$, 
the 4-manifold $\hat M$ supports a Riemannian metric of nonpositive 
curvature, and therefore is aspherical.  
There are infinitely many 4-manifolds of this form 
by recent work of M. Anderson. 

Let $\ov M$ be the compact 4-manifold with boundary 
obtained by removing disjoint horoball neighborhoods of the ideal cusp points of $M$. 
Then $\ov M$ is a strong deformation retract of $M$.  
The boundary of $\ov M$ is the disjoint union of five flat 3-tori. 
According to Ivan\v{s}i\'c \cite{Ivansic}, 
the manifold $\ov M$ is homeomorphic to the complement of the interior of a closed 
tubular neighborhood $V$ of five disjoint 2-tori $T^2_1,\ldots, T^3_5$ in $S^4$. 
The 4-manifold $V$ is the disjoint union of closed tubular neighborhoods 
$V_1, \ldots, V_5$ of $T^2_1, \ldots, T^2_5$, respectively. 
Each 4-manifold $V_i$ is homeomorphic to $D^2\times T^2$.
The boundary of $V$ is the disjoint union 
of the boundary components $T^3_1,\ldots T^3_5$ of $\ov M$ 
with $T^3_i = \partial V_i$ for each $i$. 
Let $m_i$ be a meridian of $V_i$ represented by an oriented circle for each $i$, 
and let $\kappa_i = [m_i]$ be the class of $m_i$ in $H_1(T^3_i)$ for each $i$.  
By Alexander duality, $H_1(M) \cong \integers^5$, and so $H_1(\ov M) \cong \integers^5$. 
Let 
$$\ell_i:H_1(T^3_i) \to H_1(\ov M)$$ 
be the homomorphism induced by inclusion for each $i$, 
and let $\epsilon_i = \ell_i(\kappa_i)$ for each $i$.  
By a Mayer-Vietoris sequence argument, $\epsilon_1, \ldots, \epsilon_5$   
generate $H_1(\ov M) \cong \integers^5$. 

Let $h_i:T^3_i\to T^3_i$ be an affine homeomorphism for $i=1,\ldots,5$. 
Then $h_1,\ldots,h_5$ determine an affine homeomorphism $h:\partial V \to \partial V$. 
The closed 4-manifold obtained by Dehn filling $\ov M$ according to $h_1,\ldots, h_5$ 
is the attaching space 
$$\hat M = V\cup_h\ov M.$$
Let $m_i' = h(m_i)$ for each $i$.  Then $\hat M$ is an orientable smooth 4-manifold 
whose diffeomorphism type depends only on $\pm [m_i']$ in $H_1(T^3_i)$ for each $i$. 
We also say that $\hat M$ is the closed 4-manifold obtained by Dehn surgery 
on $M$ determined by the circles $m_1',\ldots, m_5'$. 

We now prove that for infinitely Dehn surgeries on each cusp of $M$, we obtain 
a homology 4-sphere.  We use a general argument which will be simplified in \S 2 
by explicit calculations of the homomorphisms, 
$\ell_i:H_1(T^3_i) \to H_1(\ov M)$, for $i=1,\ldots,5$. 
Let $W_1 = \langle \epsilon_2,\ldots,\epsilon_5\rangle$ and let 
$$q_1: H_1(\ov M) \to H_1(\ov M)/W_1$$
be the quotient homomorphism. 
Then $q_1\ell_1: H_1(T^3_1) \to H_1(\ov M)/W_1$ is an epimorphism. 
As $H_1(\ov M)/W_1 \cong \integers$, we have a split short exact sequence
$$0 \to {\rm ker}(q_1\ell_1) \to H_1(T_1^3) \to {\rm Im}(q_1\ell_1) \to 0.$$
Hence there are generators $\kappa_{11},\kappa_{12},\kappa_{13}$ of $H_1(T_1^3)$ 
such that $\kappa_{11} =\kappa_1$ and $\kappa_{12},\kappa_{13}$ generate ${\rm ker}(q_1\ell_1)$. 

Let $b_1$ and $c_1$ be arbitrary integers. Then $\kappa_1 + b_1\kappa_{12} + c_1\kappa_{13}$ 
is a primitive element of $H_1(T^3_1)$.  Hence, there is an affine homeomorphism 
$h_1: T_1^3 \to T^3_1$ such that 
$$[h_1(m_1)] =\kappa_1 + b_1\kappa_{12} + c_1\kappa_{13}.$$
Let $m_1'=h_1(m_1)$, $\kappa_1'=[m_1']$, and $\epsilon_1'=\ell_1(\kappa_1')$. 
Then $\epsilon_1' = \epsilon_1 + \delta_1$ with $\delta_1$ in $W_1$. 
Hence $\epsilon_1',\epsilon_2,\ldots,\epsilon_5$ also generate $H_1(\ov M)$. 

Let $W_2 = \langle \epsilon_1',\epsilon_3,\ldots,\epsilon_5\rangle$ and let 
$q_2: H_1(\ov M) \to H_1(\ov M)/W_2$ be the quotient homomorphism. 
By the above argument, there are generators 
$\kappa_{21},\kappa_{22},\kappa_{23}$ of $H_1(T_2^3)$ such that 
$\kappa_{21} =\kappa_2$ and $\kappa_{22},\kappa_{23}$ generate ${\rm ker}(q_2\ell_2)$. 
Let $b_2$ and $c_2$ be arbitrary integers.  Then there is an affine homeomorphism 
$h_2: T_2^3 \to T^3_2$ such that 
$$[h_2(m_2)] = \kappa_2 + b_2\kappa_{22} + c_2\kappa_{23}.$$
Let $m_2'=h_2(m_2)$, $\kappa_2'=[m_2']$, and $\epsilon_2'=\ell_2(\kappa_2')$. 
Then $\epsilon_2' = \epsilon_2 + \delta_2$ with $\delta_2$ in $W_2$. 
Hence $\epsilon_1',\epsilon_2',\epsilon_3,\epsilon_4,\epsilon_5$ also generate $H_1(\ov M)$. 

Continuing in this way, we obtain affine homeomorphisms $h_i:T^3_i \to T^3_i$ such that 
if $m_i'=h_i(m_i)$, $\kappa_i'=[m_i']$, and $\epsilon_i'=\ell_i(\kappa_i')$, 
then $\epsilon_1',\ldots,\epsilon_5'$ generate $H_1(\ov M)$. 
Moreover, there are infinitely many choices for $h_1$, and for each choice 
of $h_1,\ldots, h_i$, there are infinitely many choices for $h_{i+1}$ for each $i=1,\ldots,4$. 
The closed 4-manifold $\hat M$ obtained by Dehn filling $\ov M$ according to $h_1,\ldots, h_5$ 
has $H_1(\hat M) = 0$ by a Mayer-Vietoris sequence argument, since 
$\epsilon_1',\ldots,\epsilon_5'$ generate $H_1(\ov M)$. 
The 4-manifold $\hat M$ is orientable, since $\ov M$ and $D^2\times T^2$ are orientable. 
By Poincar\'e duality, $H_3(\hat M) = 0$. 
By a Mayer-Vietoris sequence argument, $\chi(\hat M) = \chi(\ov M) = 2$. 
Hence we have $H_2(\hat M) = 0$. 
Therefore $\hat M$ is a homology 4-sphere for infinitely many 
Dehn surgeries on each cusp of $M$. 

As discussed by M. Anderson in \S2.1 of his paper \cite{Anderson}, 
the Gromov-Thurston $2\pi$ theorem  \cite{Bleiler}
implies that the hyperbolic metric in the interior of $\ov M$ extends 
to a Riemannian metric on $\hat M$ of nonpositive curvature 
when ${\rm length}(m_i') \geq 2\pi$ for each $i$, 
in which case, $\hat M$ is aspherical by Cartan's theorem. 
There are only finitely many homology classes of oriented circles of length less than $2\pi$ 
on the flat 3-torus $T^3_i$ for each $i$. 
Therefore $\hat M$ is an aspherical homology 4-sphere 
for infinitely many Dehn surgeries on each cusp of $M$. 
By Proposition 3.8 of Anderson \cite{Anderson}, at most finitely many 
of the 4-manifolds $\hat M$, with ${\rm length}(m_i') \geq 2\pi$ for each $i$, 
are homotopically equivalent. 
Thus there are infinitely many homotopy types of 
aspherical homology 4-spheres of the form $\hat M$.

\section{Explicit Homology Calculations} 

The most symmetric manifold on our list \cite{Ratcliffe} 
of 1171 minimum volume complete hyperbolic 4-manifolds is 
the nonorientable manifold, number 1011, 
which is denoted by $N$ in \S1. 
The hyperbolic 4-manifold $N$ has symmetry group of order 320.  
In particular, $N$ has a symmetry of order five 
that cyclically permutes its five cusps.   
Each cusp of $N$ is homeomorphic to $G\times [0,\infty)$ 
where $G$ is the first nonorientable flat 
3-manifold in Theorem 3.5.9 of Wolf \cite{Wolf}. 

Let $Q$ be the regular ideal 24-cell in the conformal ball model $B^4$ 
of hyperbolic 4-space with vertices $\pm e_i$, for $i=1,\ldots,4$ 
and 
$(\pm 1/2,\pm 1/2,\pm 1/2,\pm 1/2)$. 
The 24-cell $Q$ has 24 sides 
each of which is a regular ideal octahedra. 
The hyperbolic 4-manifold $N$ was constructed in \cite{Ratcliffe} by gluing together 
pairs of sides of the 24-cell $Q$ 
according to the side-pairing code 14FF28. 

The regular ideal 24-cell $Q$ has 24 ideal vertices. 
The side-pairing of $Q$ induces an equivalence relation 
on the ideal vertices of $Q$ whose equivalence classes are called cycles. 
The pair of vertices $\pm e_i$ form a cycle for each $i=1,\ldots, 4$, 
and the 16 vertices $(\pm 1/2,\pm 1/2,\pm 1/2,\pm 1/2)$ form a single cycle. 
The cycles of ideal vertices of $Q$ correspond to the ideal cusp points of 
the hyperbolic 4-manifold $N$.  We order the cusps of $N$ 
so that $\pm e_i$ is the cycle corresponding to the $i$th ideal cusp 
point of $N$ for $i=1,\ldots, 4$ and $(\pm 1/2,\pm 1/2,\pm 1/2,\pm 1/2)$ 
is the cycle corresponding to the fifth ideal cusp point of $N$. 

Let $\ov N$ be the compact 4-manifold with boundary obtained 
by removing disjoint horocusp neighborhoods 
of the ideal cusp points of $N$ which are invariant under 
the group of symmetries of $N$. 
The compact 4-manifold ${\ov N}$ is a strong deformation retract of $N$. 
The five components of $\partial\ov N$ are isometric 
copies of the flat 3-manifold $G$.  

Let $\ov Q$ be the truncated regular ideal 24-cell obtained 
by removing from $Q$ the disjoint horoball neighborhoods of 
the ideal vertices of $Q$ corresponding to the horocusp neighborhoods 
of the ideal cusp points of $N$ which are removed from $N$ to form ${\ov N}$.  
Then $\ov Q$ is a compact 4-dimensional polytope with 48 sides, 
24 cubical sides and 24 truncated octahedral sides. 
The side-pairing of $Q$ determines a side-pairing of the truncated 
octahedral sides of $\ov Q$ 
whose quotient space is the 4-manifold $\ov N$ with boundary. 

The side-pairing of the octahedral sides of ${\ov Q}$ 
determines a side-pairing of the 24 cubical sides of ${\ov Q}$ 
whose quotient space is the boundary of $\ov N$. 
Pairs of antipodal cubes, corresponding to $\pm e_i$, for $i=1,\ldots, 4$,  
are glued together along their sides to form the first four components of $\partial{\ov N}$, 
and the remaining 16 cubes, corresponding to the vertices $(\pm 1/2,\pm 1/2,\pm 1/2,\pm 1/2)$,  
are glued together along their sides to form the fifth component of $\partial{\ov N}$. 

The cell structure of the truncated 24-cell $\ov Q$ together 
with the side-pairing of the octahedral sides of $\ov Q$ 
determines a cell structure for $\ov N$.  
From this cell structure, we computed in \cite{Ratcliffe} 
the homology groups,  
$H_1(N)\cong \integers_2^6$, $H_2(N)\cong \integers_2^4$, and $H_3(N) = 0$. 

Ivan\v{s}i\'c's link complement \cite{Ivansic} is homeomorphic to the orientable 
double cover $M$ of $N$. 
All the symmetries of $N$ lift to symmetries of $M$. 
In particular, $M$ has a symmetry of order five 
that cyclically permutes its five cusps.  
Each cusp of $M$ is homeomorphic to $T^3\times [0,\infty)$ 
where $T^3$ is a 3-torus.  

Let $Q'$ be a regular ideal 24-cell that is obtained from $Q$ by reflecting 
in a side of $Q$. 
Then the hyperbolic 4-manifold $M$ can be constructed by gluing together 
pairs of sides of $Q$ and $Q'$. 
The side-pairing of $Q$ and $Q'$ determines a side-pairing of 
the octahedral sides of the two corresponding truncated 24-cells $\ov Q$ and ${\ov Q}\hbox{}'$ 
whose quotient space is a compact 4-manifold $\ov M$ with boundary which is 
the orientable double cover of $\ov N$. 
The compact 4-manifold ${\ov M}$ is a strong deformation retract of $M$. 
The cell structure of the two truncated 24-cells $\ov Q$ and $\ov Q\hbox{}'$ 
together with the side-pairing of the octahedral sides of $\ov Q$ and $\ov Q\hbox{}'$ 
determines a cell structure for $\ov M$.  
From this cell structure, we computed the homology groups,  
$H_1(M)\cong \integers^5$, $H_2(M)\cong \integers^{10}$, and $H_3(M)\cong \integers^4$.  
It is worth noting that this homology calculation agrees with 
calculation of the homology groups of $M$, as the complement of five disjoint 
2-tori in $S^4$, by Alexander duality.

Let $k_i:\pi_1(T_i^3) \to \pi_1(\ov M)$ be the homomorphism 
induced by inclusion for each $i=1,\ldots, 5$. 
Some care must be taken with choices of base points 
but we will suppress base points to simplify notation. 
We explicitly compute the homomorphisms $k_i:\pi_1(T_i^3) \to \pi_1(\ov M)$ 
in terms of the side-pairing transformations of the convex fundamental domain 
$Q \cup Q'$ for the hyperbolic 4-manifold $M$. 
We then explicitly compute the homomorphisms $\ell_i:H_1(T_i^3) \to H_1(\ov M)$ 
by abelianizing the homomorphisms $k_i:\pi_1(T_i^3) \to \pi_1(\ov M)$. 

We derive a group presentation for $\pi_1(\ov M)$ 
with generators the side-pairing transformations of the convex fundamental 
domain $Q \cup Q'$ for the hyperbolic 4-manifold $M$ 
by Poincar\'e's fundamental polyhedron theorem \cite{Ratcliffe0}. 
For each $i=1,\ldots, 4$,  
we form a fundamental domain for $T^3_i$ from the cubical side 
of $\ov Q$ corresponding to $e_i$ and transforms of 
three cubical sides of $\ov Q$ and $\ov Q\hbox{}'$. 
We derive a group presentation for $\pi_1(T^3_i)$ 
by Poincar\'e's fundamental polyhedron theorem. 
We then explicitly compute the homomorphism $k_i:\pi_1(T_i^3) \to \pi_1(\ov M)$ 
by rewriting the generators of the presentation for $\pi_1(T^3_i)$ 
in terms of the generators of the presentation for $\pi_1(\ov M)$ for each $i=1,\ldots,4$. 
The last homomorphism $k_5:\pi_1(T_5^3) \to \pi_1(\ov M)$ is computed in the same way 
except that we require 32 cubes to assemble a fundamental domain for $T^3_5$. 

We found generators $\kappa_{i1},\kappa_{i2},\kappa_{i3}$ for $H_1(T^3_i)$ for $i=1,\ldots,5$ 
and generators $\epsilon_1,\ldots,\epsilon_5$ for $H_1(\ov M)$ such that 
$\ell_i(\kappa_{i1}) = \epsilon_i$ for each $i$ 
and $\kappa_{i2},\kappa_{i3}$ generate ${\rm ker}(\ell_i)$ for each $i$.

Let $b_i, c_i$ be arbitrary integers for $i=1,\ldots,5$. 
By the argument at the end of \S 2, the closed 4-manifold $\hat M$ obtained 
by Dehn filling $\ov M$ according to the affine homeomorphisms $h_i:T^3_i\to T^3_i$ 
so that 
$$[h_i(m_i)] = \kappa_{i1} + b_i\kappa_{i2}+c_i\kappa_{i3},$$
for each $i$, is a homology 4-sphere; 
moreover, $\hat M$ is aspherical if the length of the circle 
$m_i' = h_i(m_i)$ is at least $2\pi$ for each $i$.

\section{The Geometry of Ivan\v{s}i\'c's link complement} 

This section will appear in version 2. 

\section{Properties of our Aspherical Homology 4-spheres} 

Smooth 4-manifolds that are homology 4-spheres have some nice 
topological properties. 
Every smooth homology 4-sphere $X$ has zero signature and 
is a spin manifold with a unique spin structure. 
By a Theorem of T. Cochran \cite{Cochran}, 
every smooth homology 4-sphere $X$ smoothly embeds in $\realnos^5$. 

Let $\hat M$ be the closed 4-manifold obtained by Dehn surgery on $M$  
determined by the circles $m'_1, \ldots, m'_5$. 
If ${\rm length}(m_i') \geq 2\pi$ for each $i$, then $\hat M$ admits 
a Riemannian metric of nonpositive curvature. 
This implies that the universal cover of $\hat M$ is diffeomorphic to $\realnos^4$ 
by Cartan's theorem; in which case $\hat M$ is aspherical, and so $\hat M$ 
is a $K(\pi,1)$.  If $\hat M$ is an aspherical homology 4-sphere, then 
$\pi_1(M)$ is an ultra super perfect group, since $H_i(\pi) = 0$ for $i=1,2,3$. 

M. Anderson \cite{Anderson} has proved that if the length of $m'_i$ 
is sufficiently large for each $i$ 
and the lengths of $m'_1,\ldots,m'_5$ are weakly balanced, in the sense that 
$${\rm max}_i({\rm length}(m_i')) \leq \exp(c\,{\rm min}_i({\rm length}(m_i'))^3)$$ 
for some small constant $c$, then the closed 4-manifold $\hat M$ 
admits an Einstein metric $g$ so that ${\rm Ric}_g = -3g$. 
Moreover the volume of $(\hat M, g)$ is less than the hyperbolic volume of $M$ 
and the volume of $(\hat M, g)$ approaches the volume of $M$ as ${\rm min}_i(m_i')$ 
goes to infinity. Thus infinitely many of the aspherical homology 4-spheres of
the form $\hat M$ admit an Einstein metric $g$ so that ${\rm Ric}_g = -3g$. 

\vspace{.2in}
{\bf Acknowledgement.} 
The authors wish to thank S. Carlip and S. Surya for helpful conversations 
about applications of hyperbolic 4-manifolds to Physics (see \cite{A-T})  
which eventually led to the discovery of our examples.


\begin{thebibliography}{99}


\bibitem{Anderson} M.T. Anderson, 
Dehn filling and Einstein metrics in higher dimensions, 
arXiv:math.DG /0303260 v3. 17 Oct 2003.

\bibitem{A-T} M. Anderson, S. Carlip, J.G. Ratcliffe, S. Surya, S.T. Tschantz, 
Peaks in the Hartle-Hawking wave function from sums over topologies, 
arXiv:gr-qc /0310002 v3. 19 Nov 2003.

\bibitem{Benedetti} R. Benedetti and C. Petronio, {\it Lectures on Hyperbolic Geometry}, 
Springer-Verlag, Berlin, 1992. 

\bibitem{Bleiler} S.A. Bleiler and C.D. Hodgson, 
Spherical space forms and Dehn filling, 
{\it Topology} {\bf 35} (1996), 809-833. 

\bibitem{Cochran} T. Cochran, 
Embedding 4-manifolds in $S^5$, 
{\it Topology} {\bf 23} (1984), 257-269. 

\bibitem{Ivansic} D. Ivan\v{s}i\'c, Hyperbolic structure on a complement of 
tori in the 4-sphere, {\it Adv. Geom.} {\bf 4} (2004).

\bibitem{Kirby} R.  Kirby, Problems in low dimensional manifold theory, 
{\it Proc. Symposia Pure Math.} {\bf 32} (1978), 273-312. 

\bibitem{Luo} F. Luo, The existence of $K(\pi,1)$ 4-manifolds which are 
rational homology 4-spheres, 
{\it Proc. Amer. Math. Soc.} {\bf 104} (1988), 1315-1321. 


\bibitem{Ratcliffe0} J.G. Ratcliffe {\it Foundations of Hyperbolic Manifolds},
Graduate Texts in Math., vol. {\bf 149}, Springer-Verlag, Berlin, Heidelberg, and
New York, 1994. 

\bibitem{Ratcliffe} J.G. Ratcliffe and S.T. Tschantz, The volume spectrum of 
hyperbolic 4-manifolds, Experiment. Math. {\bf 9}  (2000), 101-125. 

\bibitem{Thurston} W.P. Thurston, {\it The Geometry and Topology of 3-Manifolds}, 
Princeton Univ. (1979), www.msri.org/publications/books/gt3m/

\bibitem{Wolf} J.A. Wolf, {\it Spaces of Constant Curvature}, Publish
or Perish, Houston, 1974.



\end{thebibliography}
\end{document}